\documentclass[12pt,a4paper]{article}
\usepackage{amsmath,amsthm,amssymb,txfonts,bm}
\usepackage{empheq}
\usepackage[truedimen,margin=25truemm]{geometry}
\usepackage[utf8]{inputenc}

\theoremstyle{definition}
\newtheorem{definition}{Definition}
\newtheorem{thm}[definition]{Theorem}
\newtheorem{lem}[definition]{Lemma}

\newtheorem{remark}[definition]{Remark}

\def\p#1#2#3#4{{}_{2}\phi_{1}\biggl(\genfrac..{0pt}{}{{#1},\,{#2}}{#3};q,\,#4\biggr)}
\def\Q#1#2#3#4#5#6#7#8{Q\biggl(\genfrac..{0pt}{}{{#1},\,{#2}}{#3};{#4}\,;\genfrac..{0pt}{}{{#5},\,{#6}}{#7};#8\biggr)}
\def\R#1#2#3#4#5#6#7#8{R\biggl(\genfrac..{0pt}{}{{#1},\,{#2}}{#3};{#4}\,;\genfrac..{0pt}{}{{#5},\,{#6}}{#7};#8\biggr)}
\def\c#1#2#3#4#5#6#7#8{\biggl(\genfrac..{0pt}{}{{#1},\,{#2}}{#3};{#4}\,;\genfrac..{0pt}{}{{#5},\,{#6}}{#7};#8\biggr)}
\def\y#1#2#3#4#5{y_{#1}\biggl(\genfrac..{0pt}{}{{#2},\,{#3}}{#4};#5\biggr)}
\def\tQ#1#2#3#4#5#6#7#8{\tilde{Q}\biggl(\genfrac..{0pt}{}{{#1},\,{#2}}{#3};{#4}\,;\genfrac..{0pt}{}{{#5},\,{#6}}{#7};#8\biggr)}
\def\Y#1#2#3#4#5#6#7#8{Y\biggl(\genfrac..{0pt}{}{{#1},\,{#2}}{#3};{#4}\,;\genfrac..{0pt}{}{{#5},\,{#6}}{#7};#8\biggr)}

\makeatletter
\@addtoreset{equation}{section}

\makeatother

\begin{document}

\title{Basic hypergeometric identities derived from three-term relations}
\author{Yuka Yamaguchi}
\date{\today}

\maketitle

\begin{abstract}
In 2015, Ebisu presented a new method for finding hypergeometric identities based on three-term relations for the ${}_{2} F_{1}$ hypergeometric series. 
By using this method, he derived almost all of the previously known hypergeometric identities, as well as many new ones. 
In this paper, we derive several basic hypergeometric identities, including both well-known and not widely known ones, by applying a $q$-analogue of Ebisu's method to three-term relations for the ${}_{2} \phi_{1}$ basic hypergeometric series. 


\end{abstract}

\section{Introduction}
In 2015, Ebisu~\cite{Eb2} presented a new method for finding hypergeometric identities based on three-term relations for the ${}_{2} F_{1}$ hypergeometric series. 
A three-term relation is a linear relation among three ${}_{2} F_{1}$-series, in which the parameters $\alpha$, $\beta$, and $\gamma$ in each series differ by integers from those in the others. 
The coefficients in the relation are rational functions of $\alpha$, $\beta$, $\gamma$, and the variable $x$. 
Ebisu's method is as follows:  
Starting from a given three-term relation, one seeks a quadruple $(\alpha, \beta, \gamma, x)$ such that the relation reduces to a linear relation between two ${}_{2}F_{1}$-series.
This reduced relation can be regarded as a first-order difference equation with rational function coefficients.
Solving this equation leads to a hypergeometric identity.
Using this method, he derived almost all of the previously known hypergeometric identities, as well as many new ones. 
In this paper, we derive several basic hypergeometric identities, including both well-known and not widely known ones, by applying a $q$-analogue of Ebisu's method to three-term relations for the ${}_{2} \phi_{1}$ basic hypergeometric series. 

The ${}_{2} \phi_{1}$ basic hypergeometric series is defined by 
\begin{align*}
\p{a}{b}{c}{x} = {}_{2}\phi_{1} (a, b; c; q, x) := \sum_{i = 0}^{\infty} 
\frac{(a; q)_{i} (b; q)_{i}}{(q; q)_{i} (c; q)_{i}} x^{i}, 
\end{align*}
where $(a; q)_i$ denotes the $q$-shifted factorial defined by $(a; q)_{i} := (a; q)_{\infty} / (a q^{i}; q)_{\infty}$ with 
$(a; q)_{\infty} := \prod_{j = 0}^{\infty} (1 - a q^{j})$ 
and it is assumed that $\lvert q \rvert < 1$. 

Any three ${}_{2}\phi_{1}$ basic hypergeometric series, in which the parameters $a$, $b$, $c$ and the variable $x$ in each series differ multiplicatively by integer powers of $q$ from those in the others, satisfy a linear relation whose coefficients are rational functions of $a$, $b$, $c$, $q$, and $x$. 
This relation is called a three-term relation for ${}_{2}\phi_{1}$-series. 

To derive hypergeometric identities, we use three-term relations of the following form, where $k, l, m$, and $n$ are integers: 
\begin{align}\label{3tr}
\p{a q^{k}}{b q^{l}}{c q^{m}}{x q^{n}} 
= Q \cdot \p{a q}{b q}{c q}{x} + R \cdot \p{a}{b}{c}{x}, 
\end{align}
where $Q$ and $R$ are rational functions of $a, b, c, q$, and $x$. 
In \cite[Theorem~$1$]{Y1}, it is established that the pair $(Q, R)$ satisfying $(\ref{3tr})$ is uniquely determined by each $(k, l, m, n)$, and explicit expressions for $Q$ and $R$ are presented. 
The special cases with $n = 0$ are treated in \cite[Lemma~$1$, Theorem~$2$]{Y0}. 
We now describe a $q$-analogue of Ebisu's method. 
For a positive integer $N$, let 
$$
(a_N, b_N, c_N, x_N) := (aq^{k(N-1)}, bq^{l(N-1)}, cq^{m(N-1)}, xq^{n(N-1)}). 
$$ 
Then, replacing $(a,b,c,x)$ in $(\ref{3tr})$ with $(a_N, b_N, c_N, x_N)$ yields 
\begin{align}\label{3tr_N}
\p{a_{N+1}}{b_{N+1}}{c_{N+1}}{x_{N+1}} 
= Q^{(N)} \cdot \p{a_{N} q}{b_{N} q}{c_{N} q}{x_{N}} + R^{(N)} \cdot \p{a_{N}}{b_{N}}{c_{N}}{x_{N}}, 
\end{align}
where $Q^{(N)} := Q \rvert_{(a,\, b,\, c,\, x) \to (a_N,\, b_N,\, c_N,\, x_N)}$ and $R^{(N)} := R \rvert_{(a,\, b,\, c,\, x) \to (a_N,\, b_N,\, c_N,\, x_N)}$. 
Let $(a, b, c, x)$ be a solution to the system 
\begin{align}\label{Q^(N)}
Q^{(N)} = 0, \quad N = 1, 2, \ldots.
\end{align}
Then, substituting this solution into $(\ref{3tr_N})$, 
we obtain the first-order difference equation 
\begin{align*}
\p{a_{N+1}}{b_{N+1}}{c_{N+1}}{x_{N+1}} 
= R^{(N)} \cdot \p{a_{N}}{b_{N}}{c_{N}}{x_{N}}, \quad N = 1, 2, \ldots.
\end{align*}
This leads to the following equation: 
\begin{align}\label{equation}
\p{a}{b}{c}{x} 
= \frac{1}{R^{(1)} R^{(2)} \cdots R^{(N)}} \cdot \p{aq^{kN}}{bq^{lN}}{cq^{mN}}{xq^{nN}}, \quad N = 1, 2, \ldots.  
\end{align}
Manipulating this so that the right-hand side can be expressed in terms of the $q$-shifted factorial yields a basic hypergeometric identity. 
For example, if $b$ is left as a free parameter and $l > 0$, then substituting $b = q^{-lN-l'}$ with $0 \leq l' < l$ gives a closed-form expression for ${}_{2}\phi_{1}(a, q^{-lN-l'}; c; q, x)$, since the ${}_{2}\phi_{1}$-series on the right-hand side of $(\ref{equation})$ becomes a finite sum with at most $l'$ terms and can therefore be evaluated explicitly. 
Alternatively, if permissible, letting $N \to \infty$ in $(\ref{equation})$ can also lead to a basic hypergeometric identity. 

While the standard definition of the series ${}_{2}\phi_{1} (a, b; c; q, x)$ assumes that $c \neq q^{-j}$ for $j = 0, 1, \ldots$ so that the denominator factor $(c; q)_i$ does not become zero, we extend the definition to include certain exceptional cases where $a$ or $b = q^{-r}$ and $c = q^{-s}$ for integers $r$ and $s$ satisfying $0 \leq r < s$. In such cases, we define the series by
\begin{align*}
\p{a}{b}{c}{x} 
:= \sum_{i = 0}^{r} 
\frac{(a; q)_{i} (b; q)_{i}}{(q; q)_{i} (c; q)_{i}} x^{i}. 
\end{align*}

The main result of this paper is as follows. 

\begin{thm}\label{thm:1}
When $(k,l,m,n)=(0,0,0,2)$, the three-term relation $(\ref{3tr})$ leads to the identity 
\begin{align}
\p{a}{-a}{-q}{x} 
= \frac{(a^2 x;q^2)_{\infty}}{(x; q^2)_{\infty}}, 
\qquad \vert x \rvert < 1. 
\label{q-binom.2} 
\end{align}
When $(k,l,m,n)=(0,1,1,0)$, the relation $(\ref{3tr})$ yields  
\begin{align}
\p{a}{b}{c}{\frac{c}{ab}} 
= \frac{(c/a; q)_{\infty} (c/b; q)_{\infty}}{(c; q)_{\infty} (c/ab; q)_{\infty}}, 
\qquad \vert c/ab \rvert < 1. \label{q-Gauss} 
\end{align}
When $(k,l,m,n)=(0,2,2,0)$, the relation $(\ref{3tr})$ gives 
\begin{align}
\p{a}{b}{bq/a}{-\frac{q}{a}} 
= \dfrac{(-q; q)_\infty (bq;q^2)_\infty (bq^2/a^2;q^2)_\infty}{(-q/a;q)_\infty (bq/a;q)_\infty}, \qquad \lvert -q/a \rvert < 1. \label{q-Kummer}
\end{align}
When $(k,l,m,n)=(1,2,1,-1)$, the relation $(\ref{3tr})$ leads to the identities  
\begin{empheq}[left = {\displaystyle \p{a}{b}{bq/a}{-\frac{q}{a}} = \empheqlbrace \,}]{alignat = 2}
&\dfrac{(-q^{M+2};q)_N (q;q^2)_N}{(q^{M+N+2};q)_N}, 
&a=q^{M+2},\; b=q^{-2N}, \label{sv1} \\
&0, 
&a=q^{M+2},\; b=q^{-2N-1}, \label{sv2} \\
&\dfrac{(-q;q)_N (q^{M+3};q^2)_N}{q^{N(N+1)/2} (q^{M+2};q)_N}, 
&\quad a=q^{-N},\; b=q^{-M-2N-2}. \label{sv3}
\end{empheq}
When $(k,l,m,n)=(0,3,3,0)$, the relation $(\ref{3tr})$ yields 
\begin{align}
\p{\omega q}{q^{-N}}{\omega q^{-N}}{1} 
= \frac{1 - \omega^{N+1}}{1 - \omega} \frac{(q^{-N};q)_N}{(\omega q^{-N};q)_N}, 
\qquad \omega^3 = 1, \;\omega \neq 1. 
\label{sv4}
\end{align}
In the above identities, $M$ and $N$ denote arbitrary non-negative integers. 
\end{thm}

The identity~$(\ref{q-binom.2})$ is equivalent to the $q$-binomial theorem~\cite[(1.3.2)]{GR}: 
\begin{align}\label{q-binom}
\sum_{i=0}^{\infty} \frac{(a;q)_i}{(q;q)_i}x^i = \frac{(ax;q)_\infty}{(x;q)_\infty}, \qquad \lvert x \rvert < 1. 
\end{align}
The identity~$(\ref{q-Gauss})$, originally given by Heine~\cite{Heine}, is a $q$-analogue of Gauss' summation formula for ${}_{2}F_{1}(\alpha, \beta; \gamma; 1)$. 
The identity~$(\ref{q-Kummer})$, independently discovered by Bailey~\cite{Bailey} and Daum~\cite{Daum}, is a $q$-analogue of Kummer's summation formula for ${}_{2}F_{1}(\alpha, \beta; \beta + 1 - \alpha; -1)$. 
The identities~$(\ref{sv1})$--$(\ref{sv3})$ are $q$-analogues of the identities obtained by Ebisu~\cite[(1,2,2-1)-(v)]{Eb2}. 
The identity~$(\ref{sv4})$ is a special case of the identity 
\begin{align}\label{sv5}
\p{aq}{q^{-N}}{aq^{-N}}{1} 
= \frac{1-a^{N+1}}{1-a} \frac{(q^{-N};q)_N}{(aq^{-N};q)_N}, 
\end{align}
where $N$ is a non-negative integer. 
Since this identity does not appear in standard references such as \cite{GR}, we provide a brief proof in Section~$2.2$. 

We describe symmetries of the coefficient $Q$ of $(\ref{3tr})$ that helps simplify the process of finding the solutions of $(\ref{Q^(N)})$. 
Let $G$ be the group generated by the following four mappings: 
\begin{align*}
\sigma_{0} : \: &(k,l,m,n;a,b,c,x) 
\mapsto (-k,-l,-m,-n;aq^k,bq^l,cq^m,xq^n), \allowdisplaybreaks \\
\sigma_{1} : \: &(k,l,m,n;a,b,c,x) 
\mapsto (n,m-k,l+n,k;x,c/a,bx,a), \allowdisplaybreaks \\
\sigma_{2} : \: &(k,l,m,n;a,b,c,x) 
\mapsto (-k,-l,-m,k+l-m+n;q/a,q/b,q^2/c,abx/c), \allowdisplaybreaks \\
\sigma_{3} : \: &(k,l,m,n;a,b,c,x) 
\mapsto (l,k,m,n;b,a,c,x), 
\end{align*}
where the group operation is the composition of mappings. 
Let $\sigma_{4} := \sigma_{3} \sigma_{2} \sigma_{1} \sigma_{3} \sigma_{1} \sigma_{2} \sigma_{3}$, $\sigma_{5} := \sigma_{1} \sigma_{3} \sigma_{1} \sigma_{3} \sigma_{1} \sigma_{2}$, and $\sigma_{6} := \sigma_{1} \sigma_{3} \sigma_{1} \sigma_{3} \sigma_{1} \sigma_{3}$ so that they become 
\begin{align*}
\sigma_{4} : \: &(k,l,m,n;a,b,c,x) \mapsto (-n, l, m-k-n, -k;q/x,b,cq/(ax),q/a), \allowdisplaybreaks \\
\sigma_{5} : \: &(k,l,m,n;a,b,c,x) \mapsto (k-m, l-m, -m, n;aq/c, bq/c, q^2/c, x), \allowdisplaybreaks \\
\sigma_{6} : \: &(k,l,m,n;a,b,c,x) \mapsto (m-l, m-k, m, k+l-m+n;c/b, c/a, c, abx/c). 
\end{align*}
\begin{lem}\label{lem:G}
The following holds: 
\begin{align*}
G \cong \langle \sigma_{0} \rangle 
\times \left(\langle \sigma_{3}, \sigma_{4}, \sigma_{5} \rangle 
\times \langle \sigma_{6} \rangle \right) 
\cong \mathbb{Z} / 2 \mathbb{Z} 
\times \left( S_{4} \times \mathbb{Z} / 2 \mathbb{Z} \right), 
\end{align*}
where $S_{4}$ denotes the symmetric group of degree $4$. 
Thus, the order of $G$ is $2 \cdot 4! \cdot 2 = 96$. 
\end{lem}
For any $\sigma \in G$, we define its action on $Q = Q (k, l, m, n; a, b, c, x)$ by 
\begin{align*}
(\sigma Q) (k,l,m,n;a,b,c,x) := Q \left(\sigma^{-1} (k,l,m,n;a,b,c,x) \right). 
\end{align*}
Then, the following holds. 
\begin{thm}\label{thm:3}
The coefficient $Q$ of $(\ref{3tr})$ has symmetries under the action of $G$; that is, it possesses 96 symmetries. 
In fact, the following hold: 
\begin{align}
Q (k,l,m,n;a,b,c,x)
&= \frac{(-1)^{m+1} (c)_{m} (cq)_{m} (aq/c)_{k-m} (bq/c)_{l-m} (x)_{n} x^{-m}}{q^{m(m-1)/2 + mn} c^{m+n} (aq)_{k} (bq)_{l} (abxq/c)_{k+l-m+n}} \nonumber \\
&\quad \times 
(\sigma_0 Q) (k,l,m,n;a,b,c,x), \label{Q-0} \allowdisplaybreaks \\
Q (k,l,m,n;a,b,c,x) 
&= \frac{(1-b) (c-abx)}{(c-a) (c-b)} \frac{(cq)_{m-1} (xq)_{n-1}}{(aq)_{k-1} (bxq)_{l+n-1}}  
(\sigma_1 Q) (k,l,m,n;a,b,c,x), \label{Q-1} \allowdisplaybreaks \\
Q (k,l,m,n;a,b,c,x) 
&= \frac{(-1)^{m+1} (1-a) (1-b) (cq)_{m-1} (c/q^2)_{m+1} (aq/c)_{k-m} (bq/c)_{l-m}}{q^{m(m-1)/2 + (m-1)(n-1)} c^{m+n-1} (a/q)_{k+1} (b/q)_{l+1} (abxq/c)_{k+l-m+n-1}} \nonumber \\
&\quad \times (xq)_{n-1} x^{-m} (\sigma_2 Q) (k,l,m,n;a,b,c,x), \label{Q-2} \allowdisplaybreaks \\
Q (k,l,m,n;a,b,c,x) 
&= (\sigma_3 Q) (k,l,m,n;a,b,c,x). \label{Q-3}
\end{align}
\end{thm}

The symmetries~$(\ref{Q-0})$ and $(\ref{Q-3})$ are $q$-analogues of \cite[(4) and (7)]{Y2}, respectively. 
The symmetry~$(\ref{Q-2})$ is a $q$-analogue of a symmetry obtained by combining \cite[(5)--(7)]{Y2}. 
We remark that the coefficient $R$ of $(\ref{3tr})$ also has 96 symmetries (see Section~$3.3$). 

From Theorem~$\ref{thm:3}$, we find that if $(\ref{Q^(N)})$ has a solution $(a,b,c,x)$ at a certain $(k,l,m,n)$, then---except in the special cases described in Remark~$\ref{rem:4}$---it also has a solution $\sigma''(a,b,c,x)$ at $\sigma'(k,l,m,n)$, where $\sigma'$ and $\sigma''$ 
denote the restriction of $\sigma \in G$ to the $(k,l,m,n)$- and $(a,b,c,x)$-components of $(k,l,m,n;a,b,c,x)$, respectively. 
Let $G'$ be the group generated by $\sigma'_i$ ($i = 0, 1, 2, 3$). 
Then, we can take $\{(k, l, m, n) \mid 0 \leq (k+l-m)/2 \leq -n \leq k \leq l \}$ as a complete set of representatives for the quotient $\mathbb{Z}^4 / G'$ (see Section~$4$). 
The author of this paper investigated the existence of solutions to $(\ref{Q^(N)})$ for all $(k, l, m, n)$ satisfying $0 \leq (k + l - m)/2 \leq -n \leq k \leq l \leq 3$ and $m \leq 4$. 
All the resulting identities either appear in Theorem~$\ref{thm:1}$, are special cases thereof (see Section~$5$), or are trivial cases (i.e., when the ${}_{2}\phi_{1}$-series terminates after the first term, so that the sum equals $1$).

\begin{remark}\label{rem:4}
We exclude the cases
$(a,b,c,x) = (1,b,c,x)$, $(a,1,c,x)$, $(0,b,0,x)$, $(a,0,0,x)$, $(a,b,0,0)$, and $(a,b,c,0)$ from the set of solutions to $(\ref{Q^(N)})$,
since for every $\sigma \in G$, one of the following occur in these cases: 
\begin{enumerate}
\item[(i)] $(\sigma Q) (k,l,m,n;a,b,c,x) \neq 0$. 
For example, although $Q^{(N)}(0,1,1,0;1,b,c,x) = 0$ for all $N = 1, 2, \ldots$ , we have $(\sigma_2 Q)(0,1,1,0;1,b,c,x)=Q(0,-1,-1,0;q,q/b,q^2/c,bx/c)\neq 0$. 
\item[(ii)] $\sigma (k,l,m,n;a,b,c,x)$ is undefined due to division by zero. 
For instance, this occurs for $\sigma_{2} (k, l, m, n; 0, b, 0, x)$. 
\item[(iii)] Even when neither (i) nor (ii) applies, the resulting identity is either trivial, a special case of $(\ref{q-Gauss})$, or reduces to a ${}_{1}\phi_{0}$ basic hypergeometric series (namely, $(\ref{q-binom})$ or a special case thereof).
\end{enumerate}
\end{remark}

This paper is organized as follows. 
In Section~$2$, we prove Theorem~$\ref{thm:1}$ and the identity~$(\ref{sv5})$. 
Section~$3$ is devoted to the proofs of Lemma~$\ref{lem:G}$ and Theorem~$\ref{thm:3}$. 
In Section~$4$, we determine a complete set of representatives for the quotient $\mathbb{Z}^4 / G'$. 
Finally, in Section~$5$, we present our considerations and conjectures for more general $(k, l, m, n)$.

\section{Derivation of identities}
\subsection{Proof of Theorem~$\mathbf{\ref{thm:1}}$}
We prove Theorem~$\ref{thm:1}$. 
First, we derive $(\ref{q-binom.2})$. 
When $(k, l, m, n) = (0, 0, 0, 2)$, the coefficients of $(\ref{3tr})$ are given by   
\begin{align*}
Q = - \frac{(1-a)(1-b)x\{c+q-(a+b)xq \}}{(1-c)(c-abxq)}, \quad 
R = \frac{c+(1-a-b)xq}{c-abxq}
\end{align*}
according to \cite[Theorem~$1$]{Y1}. 
Since 
\begin{align*}
Q^{(N)} = - \frac{(1-a)(1-b)xq^{2(N-1)}\{c+q-(a+b)xq^{2(N-1)+1} \}}{(1-c)(c-abxq^{2(N-1)+1})}, 
\end{align*}
the quadruple $(a, b, c, x) = (a, -a, -q, x)$ satisfies $(\ref{Q^(N)})$, which leads to the equation 
\begin{align}\label{(0,0,0,2)}
\p{a}{-a}{-q}{x} 
&= \frac{1}{R^{(1)} R^{(2)} \cdots R^{(N)}} 
\cdot \p{a}{-a}{-q}{xq^{2N}} \nonumber \\
&= \frac{(a^2 x;q^2)_N}{(x;q^2)_N} \p{a}{-a}{-q}{xq^{2N}}, 
\quad N = 1, 2, \ldots. 
\end{align}
Hence, by letting $N \to \infty$, we obtain $(\ref{q-binom.2})$. 

Second, we derive $(\ref{q-Gauss})$. 
When $(k, l, m, n) = (0, 1, 1, 0)$, the coefficients of $(\ref{3tr})$ are given by   
\begin{align*}
Q = - \frac{(1-a)(c-abx)}{a-c}, \quad 
R= \frac{a(1-c)}{a-c}
\end{align*}
according to \cite[Theorem~$1$]{Y1}. 
Since the quadruple $(a, b, c, x) = (a, b, c, c/ab)$ satisfies $(\ref{Q^(N)})$, this leads to the equation 
\begin{align}\label{(0,1,1,0)}
\p{a}{b}{c}{\frac{c}{ab}} 
&= \frac{1}{R^{(1)} R^{(2)} \cdots R^{(N)}} 
\cdot \p{a}{bq^N}{cq^N}{\frac{c}{ab}} \nonumber \\
&= \frac{(c/a;q)_N}{(c;q)_N} \p{a}{bq^N}{cq^N}{\frac{c}{ab}}, 
\quad N = 1, 2, \ldots. 
\end{align}
Hence, by letting $N \to \infty$ and using $(\ref{q-binom})$, we obtain $(\ref{q-Gauss})$. 

Third, we derive $(\ref{q-Kummer})$. 
When $(k, l, m, n) = (0, 2, 2, 0)$, the coefficients of $(\ref{3tr})$ are given by   
\begin{align*}
Q &= -\frac{(1-a)(1-cq)(c-abx)\{(1-c)q+(a-bq)x\}}{(1-bq)(a-c)(a-cq)x}, \\
R &= \frac{(1-c)(1-cq)\{(1-a)cq+a(a-bq)x\}}{(1-bq)(a-c)(a-cq)x}
\end{align*}
according to \cite[Theorem~$1$]{Y1}. 
Since the quadruple $(a, b, c, x) = (a, b, bq/a, -q/a)$ satisfies $(\ref{Q^(N)})$, this leads to the equation 
\begin{align*}
\p{a}{b}{bq/a}{-\frac{q}{a}} 
&= \frac{1}{R^{(1)} R^{(2)} \cdots R^{(N)}} 
\cdot \p{a}{bq^{2N}}{bq^{2N+1}/a}{-\frac{q}{a}} \nonumber \\
&= \frac{(bq;q^2)_N (bq^2/a^2;q^2)_N }{(bq/a;q)_{2N}} 
\p{a}{bq^{2N}}{bq^{2N+1}/a}{-\frac{q}{a}}, 
\quad N = 1, 2, \ldots. 
\end{align*}
Hence, by letting $N \to \infty$ and using $(\ref{q-binom})$, we obtain $(\ref{q-Kummer})$. 

Fourth, we derive $(\ref{sv1})$--$(\ref{sv3})$. 
When $(k, l, m, n) = (1, 2, 1, -1)$, the coefficients of $(\ref{3tr})$ are given by   
\begin{align*}
Q = \frac{\{c-bq+b(1-a)x\}q}{(1-bq)(q-x)}, \quad 
R = \frac{(1-c)q}{(1-bq)(q-x)}
\end{align*}
according to \cite[Theorem~$1$]{Y1}. 
Since the quadruple $(a, b, c, x) = (a, b, bq/a, -q/a)$ satisfies $(\ref{Q^(N)})$, this leads to the equation 
\begin{align}\label{(1,2,1,-1)}
\p{a}{b}{bq/a}{-\frac{q}{a}} 
&= \frac{1}{R^{(1)} R^{(2)} \cdots R^{(N)}} 
\cdot \p{aq^N}{bq^{2N}}{bq^{N+1}/a}{-\frac{q^{1-N}}{a}} \nonumber \\
&= \frac{(-a;q)_N (bq;q^2)_N}{a^N q^{N(N-1)/2}(bq/a;q)_{N}} 
\p{aq^N}{bq^{2N}}{bq^{N+1}/a}{-\frac{q^{1-N}}{a}}, 
\quad N = 1, 2, \ldots. 
\end{align}
By substituting $b = q^{-2N}$ into $(\ref{(1,2,1,-1)})$, we obtain 
\begin{align*}
\p{a}{q^{-2N}}{q^{1-2N}/a}{-\frac{q}{a}} 
= \frac{(-a;q)_N (q^{1-2N};q^2)_N}{a^N q^{N(N-1)/2}(q^{1-2N}/a;q)_{N}} 
= \frac{(-a;q)_N (q;q^2)_N}{(aq^N;q)_N}, 
\end{align*}
which yields $(\ref{sv1})$. 
By substituting $b = q^{-2N-1}$ into $(\ref{(1,2,1,-1)})$, we obtain 
\begin{align*}
\p{a}{q^{-2N-1}}{q^{-2N}/a}{-\frac{q}{a}} 
= \frac{(-a;q)_N (q^{-2N};q^2)_N}{a^N q^{N(N-1)/2}(q^{-2N}/a;q)_{N}} 
\left\{ 1 - \frac{(1-aq^N)(1-q^{-1})}{(1-q)(1-q^{-N}/a)} \frac{q^{1-N}}{a} \right\}
=0, 
\end{align*}
which yields $(\ref{sv2})$. 
By substituting $a = q^{-N}$ into $(\ref{(1,2,1,-1)})$, we obtain 
\begin{align*}
\p{q^{-N}}{b}{bq^{N+1}}{-q^{N+1}} 
= \frac{(-q^{-N};q)_N (bq;q^2)_N}{q^{-N(N+1)/2}(bq^{N+1};q)_{N}} 
= \frac{(-q;q)_N (q^{1-2N}/b;q^2)_N}{q^{N(N+1)/2} (q^{-2N}/b;q)_N}, 
\end{align*}
which yields $(\ref{sv3})$. 
We note that $(\ref{sv1})$--$(\ref{sv3})$ can also be derived from $(\ref{q-Kummer})$ through an appropriate reformulation. 

Finally, we derive $(\ref{sv4})$. 
When $(k, l, m, n) = (0,3,3,0)$, the coefficients of $(\ref{3tr})$ are given by   
\begin{align*}
Q &= -\frac{(1-a)(c-abx)(cq;q)_2}{a^3 (bq;q)_2 (c/a;q)_3 x^2} \\
&\quad \times \{(1-c)(1-cq)q^2+(1-c)(a-bq^2)xq -(1-a)(b-c)xq^2+(a-bq)(a-bq^2)x^2 \}, \\
R &= \frac{(c;q)_3}{a^3 (bq;q)_2 (c/a;q)_3 x^2} \\
&\quad \times \{(1-a)(1-cq)cq^2+c(1-a)(a-bq^2)xq -a(1-a)(b-c)xq^2+a(a-bq)(a-bq^2)x^2 \}
\end{align*}
according to \cite[Theorem~$1$]{Y1}. 
Since the quadruple $(a, b, c, x) = (\omega q, b, \omega b, 1)$ satisfies $(\ref{Q^(N)})$, this leads to the equation 
\begin{align}\label{(0,3,3,0)}
\p{\omega q}{b}{\omega b}{1} 
&= \frac{1}{R^{(1)} R^{(2)} \cdots R^{(N)}} 
\cdot \p{\omega q}{bq^{3N}}{\omega bq^{3N}}{1} \nonumber \\
&= \frac{(b;q)_{3N}}{(\omega b;q)_{3N}} 
\p{\omega q}{bq^{3N}}{\omega bq^{3N}}{1}, 
\quad N = 1, 2, \ldots, 
\end{align}
where it is assumed that $b = q^{-3N-j}$ with $j=0, 1, \ldots$ to ensure convergence. 
By substituting $b = q^{-3N}$, $q^{-3N-1}$, and $q^{-3N-2}$ into $(\ref{(0,3,3,0)})$, we obtain $(\ref{sv4})$.

\subsection{Generalization of the last identity in Theorem~$\mathbf{\ref{thm:1}}$}
We provide a brief proof of $(\ref{sv5})$. 
Suppose that $\lvert x \rvert < 1$ and $\lvert ax \rvert <1$. 
Then, it holds that 
\begin{align}\label{pf:sv5-1}
\frac{1}{1-x} \frac{1}{1-ax} 
= \sum_{i=0}^{\infty}x^i \sum_{j=0}^{\infty}(ax)^j
= \sum_{N=0}^{\infty}\sum_{i=0}^{N}a^{N-i}x^N 
= \sum_{N=0}^{\infty} \frac{1-a^{N+1}}{1-a}x^N. 
\end{align}
On the other hand, from $(\ref{q-binom})$, it follows that 
\begin{align}\label{pf:sv5-2}
\frac{1}{1-x} \frac{1}{1-ax} 
&= \frac{(xq;q)_\infty}{(x;q)_\infty} \frac{(axq;q)_\infty}{(ax;q)_\infty} 
= \frac{(axq;q)_\infty}{(x;q)_\infty} \frac{(xq;q)_\infty}{(ax;q)_\infty} 
= \sum_{i=0}^{\infty} \frac{(aq;q)_i}{(q;q)_i}x^i 
\sum_{j=0}^{\infty} \frac{(q/a;q)_j}{(q;q)_j}(ax)^j \nonumber \\
&= \sum_{N=0}^{\infty} \sum_{i=0}^{N} \frac{(aq;q)_i (q/a;q)_{N-i}}{(q;q)_i (q;q)_{N-i}} a^{N-i} x^N 
= \sum_{N=0}^{\infty} \sum_{i=0}^{N} \frac{(aq;q)_i (q^{-N};q)_{i} (aq^{-N};q)_N}{(q;q)_i (aq^{-N};q)_{i} (q^{-N};q)_N} x^N.  
\end{align}
By comparing the coefficients of $x^N$ in $(\ref{pf:sv5-1})$ and $(\ref{pf:sv5-2})$, we obtain $(\ref{sv5})$. 

\section{Symmetries of the coefficients}
\subsection{Proof of Lemma~$\mathbf{\ref{lem:G}}$}
We prove Lemma~$\ref{lem:G}$. 
Let $H$ and $K$ be the subgroups of $G$ defined by $H := \langle \sigma_{1}, \sigma_{2}, \sigma_{3} \rangle$ and $K := \langle \sigma_{3}, \sigma_{4}, \sigma_{5} \rangle$. Let ${\rm id}_{G}$ denote the identity element of $G$. First, from $\sigma_0 \sigma_i = \sigma_i \sigma_0$ ($i = 1, 2, 3$) and $\sigma_{0}^{2} = {\rm id}_{G}$, we have $G \cong \langle \sigma_{0} \rangle \times H \cong \mathbb{Z} / 2 \mathbb{Z} \times H$. Next, from $\sigma_1 = \sigma_3 \sigma_4 \sigma_5 \sigma_4 \sigma_3 \sigma_6$ and $\sigma_2 = \sigma_3 \sigma_5 \sigma_6$, it holds that $H = \langle \sigma_{3}, \sigma_{4}, \sigma_{5}, \sigma_{6} \rangle$. 
In addition, from $\sigma_6 \sigma_i = \sigma_i \sigma_6$ ($i = 3, 4, 5$) and $\sigma_{6}^{2} = {\rm id}_{G}$, we obtain $H \cong K \times \mathbb{Z} / 2 \mathbb{Z}$. 
Finally, we prove $K \cong S_4$. Let $s_i \in S_4$ ($i = 1, 2, 3$) be the adjacent transposition $(i\;\; i+1)$. Then, the following holds: 
\begin{align*}
S_4 = \langle s_1, s_2, s_3 \mid s_1^2 = s_2^2 = s_3^2 = (),\; s_1 s_3 = s_3 s_1,\; (s_1 s_2)^3 = (s_2 s_3)^3 = () \rangle, 
\end{align*}
where $()$ denotes the identity element of $S_4$. 
On the other hand, $\sigma_i$ ($i = 3, 4, 5$) have the following relations: 
\begin{align*}
\sigma_3^2 = \sigma_4^2 = \sigma_5^2 = {\rm id}_{G}, \quad \sigma_3 \sigma_5 = \sigma_5 \sigma_3,\quad (\sigma_3 \sigma_4)^3 = (\sigma_4 \sigma_5)^3 = {\rm id}_{G}. 
\end{align*}
Therefore, we obtain an isomorphism $K \to S_4$ defined by $\sigma_{i+2} \mapsto s_i$ ($i = 1, 2, 3$). 
The lemma is proved. 

\begin{remark}
The group $G$ in this paper and the group $G$ in \cite{Y2} are isomorphic. 
In order to avoid confusion, we write $\sigma_i$ ($0 \leq i \leq 5$) appearing in \cite{Y2} as $\tau_i$. 
In \cite[Lemma~$2$]{Y2}, it is shown that 
\begin{align*}
G := \langle \tau_0, \tau_1, \tau_2, \tau_3 \rangle 
\cong \mathbb{Z} / 2 \mathbb{Z} \times \langle \tau_1, \tau_2, \tau_3 \rangle 
\end{align*}
with $\langle \tau_1, \tau_2, \tau_3 \rangle \cong (S_{3} \ltimes (\mathbb{Z} / 2 \mathbb{Z})^{3})$. 
We show that $\langle \tau_1, \tau_2, \tau_3 \rangle \cong S_4 \times \mathbb{Z} / 2 \mathbb{Z}$. 
From $\tau_3 = \tau_2 \tau_4 \tau_2 \tau_4 \tau_5$, $\tau_5 \tau_i = \tau_i \tau_5$ ($i = 1, 2, 4$), and $\tau_5^2 = {\rm id}_G$, it holds that $\langle \tau_1, \tau_2, \tau_3 \rangle = \langle \tau_1, \tau_2, \tau_4, \tau_5 \rangle \cong \langle \tau_1, \tau_2, \tau_4 \rangle \times \mathbb{Z} / 2 \mathbb{Z}$. 
Moreover, the homomorphism $\langle \tau_1, \tau_2, \tau_4 \rangle \to S_4$ defined by $\tau_1 \mapsto s_1$, $\tau_2 \mapsto s_2$, and $\tau_4 \mapsto s_1 s_3$ is an isomorphism. 
\end{remark}

\subsection{Proof of Theorem~$\mathbf{\ref{thm:3}}$}
We prove Theorem~$\ref{thm:3}$. 
In this section, we write the coefficients of $(\ref{3tr})$ as 
\begin{align*}
Q = Q \c{k}{l}{m}{n}{a}{b}{c}{x}, \quad 
R = R \c{k}{l}{m}{n}{a}{b}{c}{x}. 
\end{align*}

First, since ${}_{2}\phi_{1}(a, b; c; q, x)$ is symmetric with respect to the exchange of $a$ and $b$, it follows that $(\ref{Q-3})$ is true. 

Next, we prove $(\ref{Q-1})$. 
By replacing $(k,l,m,n;a,b,c,x)$ in $(\ref{3tr})$ with $(n,\allowbreak m-k,\allowbreak l+n,\allowbreak k;\allowbreak x,\allowbreak c/a,\allowbreak bx,a)$, we have  
\begin{align}\label{3tr'}
\p{xq^n}{cq^{m-k}/a}{bxq^{l+n}}{aq^k} 
&= \Q{n}{m-k}{l+n}{k}{x}{c/a}{bx}{a}\, \p{xq}{cq/a}{bxq}{a} \nonumber \\
&\quad+ \R{n}{m-k}{l+n}{k}{x}{c/a}{bx}{a}\, \p{x}{c/a}{bx}{a}. 
\end{align}
Applying Heine's \cite{Heine} (also available in \cite[(1.4.1)]{GR}) transformation formula
\begin{align*}
\p{a}{b}{c}{x} = \frac{(a)_\infty (bx)_\infty}{(c)_\infty (x)_\infty} \p{x}{c/a}{bx}{a}
\end{align*}
to the three ${}_{2}\phi_{1}$'s in $(\ref{3tr'})$ yields 
\begin{align*}
\frac{(xq^n)_\infty (cq^m)_\infty}{(bxq^{l+n})_\infty (aq^k)_\infty} \p{a q^{k}}{b q^{l}}{c q^{m}}{x q^{n}} 
&= \Q{n}{m-k}{l+n}{k}{x}{c/a}{bx}{a} \frac{(xq)_\infty (cq)_\infty}{(bxq)_\infty (a)_\infty} \p{a}{b}{cq}{xq} \\
&\quad + \R{n}{m-k}{l+n}{k}{x}{c/a}{bx}{a} \frac{(x)_\infty (c)_\infty}{(bx)_\infty (a)_\infty} \p{a}{b}{c}{x}. 
\end{align*}
Moreover, replacing ${}_{2}\phi_{1}(a, b; cq; q, xq)$ on the right-hand side of this by using the three-term relation 
\begin{align*}
\p{a}{b}{cq}{xq} 
&= \frac{(1-a)(1-b)(c-abx)}{(c-a)(c-b)} \p{aq}{bq}{cq}{x} \\
&\quad + \left\{1 - \frac{(1-a)(1-b)cq}{(c-a)(c-b)} \right\} \p{a}{b}{c}{x}
\end{align*}
and multiplying both sides of the resulting one by $(bxq^{l+n})_\infty (aq^k)_\infty / \{ (xq^n)_\infty (cq^m)_\infty \}$ give the following three-term relation: 
\begin{align*}
&\p{a q^{k}}{b q^{l}}{c q^{m}}{x q^{n}} \\
&= \Q{n}{m-k}{l+n}{k}{x}{c/a}{bx}{a} \frac{(xq)_{n-1} (cq)_{m-1}}{(bxq)_{l+n-1} (aq)_{k-1}} \frac{(1-b)(c-abx)}{(c-a)(c-b)} \p{aq}{bq}{cq}{x} \\
&\quad + R' \cdot \p{a}{b}{c}{x}, 
\end{align*}
where $R'$ is a certain rational function of $a, b, c, q$, and $x$. 
Hence, from the uniqueness of the coefficients of $(\ref{3tr})$, we obtain $(\ref{Q-1})$. 

To prove $(\ref{Q-0})$ and $(\ref{Q-2})$, we use the expression for $Q$ provided in \cite[(3.2) and (3.9)]{Y1}: 
\begin{align}
Q &= \frac{(c q)_{m - 1}}{(a q)_{k - 1} (b q)_{l - 1}} \tilde{Q}, \label{Q=tQ}\\
\tilde{Q} &= \tQ{k}{l}{m}{n}{a}{b}{c}{x}
= \Y{k}{l}{m}{n}{a}{b}{c}{x} \,\Bigg/\, \Y{1}{1}{1}{0}{a}{b}{c}{x}, \label{tQ}
\end{align}
where $Y$ is the infinite series defined by 
\begin{align*}
\Y{k}{l}{m}{n}{a}{b}{c}{x} &:= 
\y{1}{a q^{k}}{b q^{l}}{c q^{m}}{x q^{n}} 
\y{2}{a}{b}{c}{x} 
- \y{2}{a q^{k}}{b q^{l}}{c q^{m}}{x q^{n}} 
\y{1}{a}{b}{c}{x}, \\
\y{1}{a}{b}{c}{x} &:= \frac{(q)_\infty (c)_\infty}{(a)_\infty (b)_\infty} \p{a}{b}{c}{x}, \\
\y{2}{a}{b}{c}{x} &:= \frac{(q)_\infty (q^2/c)_\infty}{(aq/c)_\infty (bq/c)_\infty} x^{1 - \gamma} \p{aq/c}{bq/c}{q^2/c}{x}
\end{align*}
with $c = q^\gamma$. Concerning the denominator of $(\ref{tQ})$, the following holds \cite[(3.20)]{Y1}: 
\begin{align}\label{Y(1110)}
\Y{1}{1}{1}{0}{a}{b}{c}{x} 
= - \frac{(q)_{\infty}^{2} (c)_{\infty} (q/c)_{\infty} (abxq/c)_{\infty} x^{- \gamma}}{(a)_{\infty} (b)_{\infty} (aq/c)_{\infty} (bq/c)_{\infty} (x)_{\infty}}. 
\end{align}

We prove $(\ref{Q-0})$. From the definition of $Y$, it holds that 
\begin{align*}
\Y{-k}{-l}{-m}{-n}{aq^k}{bq^l}{cq^m}{xq^n} = - \Y{k}{l}{m}{n}{a}{b}{c}{x}, 
\end{align*}
so it follows from $(\ref{tQ})$ that 
\begin{align*}
\tQ{-k}{-l}{-m}{-n}{aq^k}{bq^l}{cq^m}{xq^n} 
= - \frac{\displaystyle \Y{1}{1}{1}{0}{a}{b}{c}{x}}{\displaystyle \Y{1}{1}{1}{0}{aq^k}{bq^l}{cq^m}{xq^n}}
\tQ{k}{l}{m}{n}{a}{b}{c}{x}. 
\end{align*}
By rewriting both sides by using $(\ref{Q=tQ})$ and $(\ref{Y(1110)})$, we obtain $(\ref{Q-0})$. 

Finally, we prove $(\ref{Q-2})$. From Heine's \cite{Heine} (also available in \cite[(1.4.3)]{GR}) transformation formula 
\begin{align*}
\p{a}{b}{c}{x} = \frac{(abx/c)_\infty}{(x)_\infty} \p{c/a}{c/b}{c}{\frac{abx}{c}}, 
\end{align*}
it holds that 
\begin{align*}
&\y{1}{q^{1-k}/a}{q^{1-l}/b}{q^{2-m}/c}{\frac{abxq^{k+l-m+n}}{c}} \\
&= \frac{(aq^{k-m+1}/c)_\infty (bq^{l-m+1}/c)_\infty (xq^n)_\infty}{(q^{1-k}/a)_\infty (q^{1-l}/b)_\infty (abxq^{k+l-m+n}/c)_\infty} (xq^n)^{\gamma+m-1} 
\y{2}{a q^{k}}{b q^{l}}{c q^{m}}{x q^{n}}, \\
&\y{2}{q^{1-k}/a}{q^{1-l}/b}{q^{2-m}/c}{\frac{abxq^{k+l-m+n}}{c}} \\
&= \frac{(aq^k)_\infty (bq^l)_\infty (xq^n)_\infty}{(cq^{m-k}/a)_\infty (cq^{m-l}/b)_\infty (abxq^{k+l-m+n}/c)_\infty} \left(\frac{abxq^{k+l-m+n}}{c} \right)^{\gamma+m-1} 
\y{1}{a q^{k}}{b q^{l}}{c q^{m}}{x q^{n}}. 
\end{align*}
Thus, we have 
\begin{align*}
\Y{-k}{-l}{-m}{k+l-m+n}{q/a}{q/b}{q^2/c}{\frac{abx}{c}} 
= - \lambda \cdot \Y{k}{l}{m}{n}{a}{b}{c}{x}, 
\end{align*}
where
\begin{align*}
\lambda 
:= \frac{(aq^{k-m+1}/c)_\infty (bq^{l-m+1}/c)_\infty (xq^n)_\infty (a)_\infty (b)_\infty (x)_\infty}{(q^{1-k}/a)_\infty (q^{1-l}/b)_\infty (abxq^{k+l-m+n}/c)_\infty (c/a)_\infty (c/b)_\infty (abx/c)_\infty} (xq^n)^{\gamma+m-1} \left(\frac{abx}{c} \right)^{\gamma-1}. 
\end{align*}
From this and $(\ref{tQ})$, we obtain 
\begin{align*}
&\tQ{-k}{-l}{-m}{k+l-m+n}{q/a}{q/b}{q^2/c}{\frac{abx}{c}} \\
&\quad = - \lambda \, \frac{\displaystyle \Y{1}{1}{1}{0}{a}{b}{c}{x}}{\displaystyle \Y{1}{1}{1}{0}{q/a}{q/b}{q^2/c}{\frac{abx}{c}}} \tQ{k}{l}{m}{n}{a}{b}{c}{x}. 
\end{align*}
\begin{align*}
\tQ{-k}{-l}{-m}{k+l-m+n}{q/a}{q/b}{q^2/c}{\frac{abx}{c}} 
= - \lambda \, \frac{\displaystyle \Y{1}{1}{1}{0}{a}{b}{c}{x}}{\displaystyle \Y{1}{1}{1}{0}{q/a}{q/b}{q^2/c}{\frac{abx}{c}}} \tQ{k}{l}{m}{n}{a}{b}{c}{x}. 
\end{align*}
By rewriting both sides by using $(\ref{Q=tQ})$ and $(\ref{Y(1110)})$, we obtain $(\ref{Q-2})$. 

\subsection{Symmetries of another coefficient: additional remarks}
Since the coefficients $Q$ and $R$ of $(\ref{3tr})$ have the relation \cite[Corollary~$2$]{Y1} 
\begin{align}\label{Q=R}
Q (k-1, l-1, m-1, n; aq, bq, cq, x)
= \frac{(1-aq) (1-bq) x (c-abxq)}{(1-c) (1-cq)} 
R (k, l, m, n; a, b, c, x), 
\end{align}
we find that $R$ also has $96$ symmetries. 
In fact, from Theorem~$\ref{thm:3}$ and $(\ref{Q=R})$, we obtain the following: 
\begin{align}
R (k, l, m, n; a, b, c, x)
&= \frac{(-1)^{m} (c)_{m-1} (cq)_{m-1} (aq/c)_{k-m} (bq/c)_{l-m} (x)_{n} x^{1-m}}{q^{m(m-1)/2 + (m-1)(n-1)} c^{m+n-1} (aq)_{k-1} (bq)_{l-1} (abxq/c)_{k+l-m+n-1}} \nonumber \\
&\quad \times 
(\tau \sigma_0 \tau^{-1} R) (k, l, m, n; a, b, c, x), \label{R-0} \allowdisplaybreaks \\
R (k, l, m, n; a, b, c, x)
&= \frac{(c)_{m} (x)_{n}}{(aq)_{k-1} (bx)_{l+n}}  
(\tau \sigma_1 \tau^{-1} R) (k, l, m, n; a, b, c, x), \label{R-1} \allowdisplaybreaks \\
R (k, l, m, n; a, b, c, x)
&= \frac{(-1)^{m} (c)_{m-1} (cq)_{m-1} (aq/c)_{k-m} (bq/c)_{l-m} (x)_{n} x^{1-m}}{q^{m(m-1)/2 + (m-1)(n-1)} c^{m+n-1} (aq)_{k-1} (bq)_{l-1} (abxq/c)_{k+l-m+n-1}} \nonumber \\
&\quad \times 
(\tau \sigma_2 \tau^{-1} R) (k, l, m, n; a, b, c, x), \label{R-2} \allowdisplaybreaks \\
R (k, l, m, n; a, b, c, x)
&= (\tau \sigma_3 \tau^{-1} R) (k, l, m, n; a, b, c, x), \label{R-3}
\end{align}
where 
\begin{align*}
\tau : \: &\c{k}{l}{m}{n}{a}{b}{c}{x} \mapsto \c{k+1}{l+1}{m+1}{n}{a/q}{b/q}{c/q}{x}. 
\end{align*}
This means that $R$ has symmetries under the action of the group generated by $\tau \sigma_i \tau^{-1}$ ($i = 0, 1, 2, 3$), whose order is 96. 

\section{Complete set of representatives for $\mathbf{\mathbb{Z}^4 / G'}$}
We determine a complete set of representatives for the quotient $\mathbb{Z}^4 / G'$.

Let $\sigma'$ denote the restriction of $\sigma \in G$ to the $(k,l,m,n)$-components of $(k,l,m,n;a,b,c,x)$, and let $G'$ be the group generated by $\sigma'_i$ ($i = 0, 1, 2, 3$). 
From $\sigma'_0 \sigma'_i = \sigma'_i \sigma'_0$ ($i=1,2,3$), $\sigma'_1 = \sigma'_3 \sigma'_4 \sigma'_5 \sigma'_4 \sigma'_3 \sigma'_6$, and $\sigma'_2 = \sigma'_3 \sigma'_5 \sigma'_6$, 
it follows that $G' := \langle \sigma'_0, \sigma'_1, \sigma'_2, \sigma'_3 \rangle = \langle \sigma'_0, \sigma'_3, \sigma'_4, \sigma'_5, \sigma'_6 \rangle$.  
Let $T$ be the coordinate transformation on $\mathbb{Z}^{4}$ defined by $T :(k,l,m,n) \mapsto (\lambda_1, \lambda_2, \lambda_3, \lambda_4) := (k,l,-n,k+l-m)$. 
Then, from
\begin{align*}
\sigma'_{0}(k,l,m,n) &= (-k,-l,-m,-n), \allowdisplaybreaks \\
\sigma'_{3}(k,l,m,n) &= (l,k,m,n), \allowdisplaybreaks \\
\sigma'_{4}(k,l,m,n) &= (-n,l,m-k-n,-k), \allowdisplaybreaks \\
\sigma'_{5}(k,l,m,n) &= (k-m,l-m,-m,n), \allowdisplaybreaks \\
\sigma'_{6}(k,l,m,n) &= (m-l,m-k,m,k+l-m+n), 
\end{align*}
we have
\begin{align*}
T \sigma'_{0} T^{-1}(\lambda_1, \lambda_2, \lambda_3, \lambda_4) &= (-\lambda_1, -\lambda_2, -\lambda_3, -\lambda_4), \allowdisplaybreaks \\
T \sigma'_{3} T^{-1}(\lambda_1, \lambda_2, \lambda_3, \lambda_4) &= (\lambda_2, \lambda_1, \lambda_3, \lambda_4), \allowdisplaybreaks \\
T \sigma'_{4} T^{-1}(\lambda_1, \lambda_2, \lambda_3, \lambda_4) &= (\lambda_3, \lambda_2, \lambda_1, \lambda_4), \allowdisplaybreaks \\
T \sigma'_{5} T^{-1}(\lambda_1, \lambda_2, \lambda_3, \lambda_4) &= (\overline{\lambda}_2, \overline{\lambda}_1, \lambda_3, \lambda_4), \allowdisplaybreaks \\
T \sigma'_{0} \sigma'_{6} T^{-1}(\lambda_1, \lambda_2, \lambda_3, \lambda_4) &= (\overline{\lambda}_1, \overline{\lambda}_2, \overline{\lambda}_3, \lambda_4)
\end{align*}
with $\overline{\lambda}_i := \lambda_4 - \lambda_i$ ($i = 1, 2, 3$). 
Therefore, 
\begin{align*}
\{(\lambda_1, \lambda_2, \lambda_3, \lambda_4) \mid \lambda_4 \geq 0, \, \overline{\lambda}_3 \leq \lambda_3 \leq \lambda_1 \leq \lambda_2 \}
\end{align*}
is a complete set of representatives for the quotient of the action of $T G' T^{-1}$ on $\mathbb{Z}^4$. 
Noting that the condition $\overline{\lambda}_3 \leq \lambda_3$ is equivalent to $\lambda_4/2 \leq \lambda_3$, it follows that 
\begin{align*}
\{(k, l, m, n) \mid 0 \leq (k+l-m)/2 \leq -n \leq k \leq l \}
\end{align*}
is a complete set of representatives for the quotient $\mathbb{Z}^4 / G'$.

\section{Considerations and conjectures for more general cases}
We present our considerations and conjectures concerning more general $(k,l,m,n)$. 

When $(k,l,m,n)=(l,l,0,n)$, with $n$ an even integer, it follows from \cite[Lemma~$15$]{Y1} that the quadruple $(a, b, c, x) = (a, -a, -q, x)$ satisfies $(\ref{Q^(N)})$. 
This leads to the equation 
\begin{align*}
\p{a}{-a}{-q}{x} 
&= \frac{1}{R^{(1)} R^{(2)} \cdots R^{(N)}} 
\cdot \p{aq^{lN}}{-aq^{lN}}{-q}{xq^{nN}}, 
\quad N = 1, 2, \ldots. 
\end{align*}
Therefore, any three-term relation $(\ref{3tr})$ in the case $(k,l,m,n)=(l,l,0,n)$, where $l$ is a non-negative integer and $n$ is a positive even integer, would lead to the identity~$(\ref{q-binom.2})$, as demonstrated above in the special case $(k,l,m,n)=(0,0,0,2)$. 

When $k+l-m+n=0$, the rational function $Q^{(N)}$ has the factor $(c-abx)$, according to \cite[Theorem~$1$]{Y1}, so the quadruple $(a, b, c, x) = (a, b, c, c/ab)$ satisfies $(\ref{Q^(N)})$. 
This leads to the equation 
\begin{align*}
\p{a}{b}{c}{\frac{c}{ab}} 
&= \frac{1}{R^{(1)} R^{(2)} \cdots R^{(N)}} 
\cdot \p{aq^{kN}}{bq^{lN}}{cq^{mN}}{\frac{cq^{nN}}{ab}}, 
\quad N = 1, 2, \ldots. 
\end{align*}
Therefore, any three-term relation $(\ref{3tr})$ in the case $k+l-m+n=0$, with $k,l,n \geq 0$ and $m>0$, would lead to the identity~$(\ref{q-Gauss})$, as demonstrated above in the special case $(k,l,m,n)=(0,1,1,0)$. 

When $(k,l,m,n)=(k,l,l-k,-k)$, with $l$ an even integer, 
we can verify that the quadruple $(a, b, c, x) = (a, b, bq/a, -q/a)$ satisfies $(\ref{Q^(N)})$ by using \cite[Lemma~$15$]{Y1} together with the $q$-Kummer summation formula~$(\ref{q-Kummer})$. 
This leads to the equation
\begin{align*}
\p{a}{b}{bq/a}{-\frac{q}{a}} 
&= \frac{1}{R^{(1)} R^{(2)} \cdots R^{(N)}} 
\cdot \p{aq^{kN}}{bq^{lN}}{bq^{(l-k)N+1}/a}{-\frac{q^{1-kN}}{a}}, 
\quad N = 1, 2, \ldots. 
\end{align*}
Therefore, as shown above in the case $(k,l,m,n)=(1,2,1,-1)$, any three-term relation $(\ref{3tr})$ in the case $(k,l,m,n)=(k,l,l-k,-k)$, with $l$ a positive even integer, would lead to the identities~$(\ref{sv1})$ and~$(\ref{sv2})$, or to special cases thereof. 
In particular, when $k>0$, it would also lead to the identity~$(\ref{sv3})$, or to a special case thereof. 

We also conjecture that, when $(k,l,m,n)=(0,l,l,0)$ with $l \geq 4$, 
the three-term relation $(\ref{3tr})$ leads to a special case of the identity~$(\ref{sv5})$ with $a=\zeta_{l}$, a primitive $l$-th root of unity; the case $l=3$ was demonstrated above.

\clearpage

\section*{Acknowledgment}
This work was supported by JSPS KAKENHI Grant Number JP25KJ0266. 

\section*{Data Availability Statement}
This paper has no associated data. 

\section*{Competing interests}
The author has no competing interests to declare that are relevant to the content of this paper. 

\bibliography{reference}

\medskip
\begin{flushleft}
Faculty of Education \\ 
University of Miyazaki \\ 
1-1 Gakuen Kibanadai-nishi \\ 
Miyazaki 889-2192 Japan \\ 
{\it Email address}: y-yamaguchi@miyazaki-u.ac.jp 
\end{flushleft}

\end{document}